\documentclass[a4paper]{amsart}
\usepackage[all]{xy}
\usepackage{setspace,amssymb}
\newtheorem{thm}{Theorem}[section]
\newtheorem{lem}[thm]{Lemma}
\newtheorem{prop}[thm]{Proposition}
\newtheorem{cor}[thm]{Corollary}
\newtheorem*{MT}{Main Theorem}
\newtheorem*{IC}{Important Corollary}
\theoremstyle{definition}

\theoremstyle{remark}
\newtheorem{rmk}[thm]{Remark}
\newtheorem{ex}[thm]{Example}
\DeclareMathAlphabet{\mathbfit}{OT1}{cmr}{bx}{it}

\newcommand{\Z}{\mathbf{Z}}
\newcommand{\Q}{\mathbf{Q}}

\newcommand{\F}[1]{\mathbf{F}_{\!{#1}}}
\newcommand{\Zmod}[1]{\Z/{#1}}

\newcommand{\op}[1]{\mathcal{{#1}}}
\newcommand{\free}[2]{\mathbf{{#1}}({#2})}

\newcommand{\Sq}[1]{\mathrm{Sq}^{{#1}}}
\begin{document}

\title[$E_{\infty}$ Cell Models]{$E_{\infty}$ Cell Models
for Free
and Based Loop Space Cohomology}
\author{David Chataur}
\address{Centre de Recerca Matem\`atica \\
    Apartat 50 \\
    E-08193 Bellaterra \\
    Spain}
\email{DChataur@crm.es}
\author{Jonathan A. Scott}
\address{Laboratoire J.-A. Dieudonn\'e \\
    Universit\'e de Nice Sophia-Antopolis \\
    Parc Valrose \\
    06108 Nice Cedex \\
    France}
\curraddr{UFR Math\'{e}matiques Pures et Appliqu\'{e}es \\
    Universit\'{e} des Sciences et Technologies de Lille \\
    59655 Villeneuve d'Ascq CEDEX \\
    France}
\email{Jonathan.Scott@agat.univ-lille1.fr}
\date{\today}
\thanks{The second author was supported in part by an NSERC Post-Doctoral
Fellowship and a European TMR Post-Doctoral Fellowship.}

\begin{abstract}
We construct $E_{\infty}$
cell algebra models for the free and based loop spaces on a
simply-connected topological space.
Techniques from rational homotopy theory are exploited throughout.
\end{abstract}

\maketitle

\section{Introduction}\label{sec:intro}

The purpose of this paper is to use methods from rational homotopy theory
to construct $E_{\infty}$ \emph{cell models} for the normalised singular
cochain algebras of the free and based loop spaces on a simply-connected
topological space.
The normalised singular cochains with coefficients in a commutative ring $R$,
denoted $N^{*}(-)$, is a functor from spaces to $E_{\infty}$ algebras
over $R$ (see, for example,~\cite{berger-fresse:01}).
The $E_{\infty}$ algebras over $R$ form a closed model category in which
the fibrations are the surjections and the weak equivalences are the
quasi-isomorphisms.
The cofibrations are retracts of
\emph{cell extensions} (see Section~\ref{sec:operads}), which are the
$E_{\infty}$ analogues of the KS extensions, or relative Sullivan algebras,
of rational homotopy theory.
A \emph{cell model} for a morphism
$\varphi:A\rightarrow B$ is a factorisation of $\varphi$ as a cell extension
followed by a weak equivalence.
A \emph{cell model} for $A$ is a cell model of the unit map.
Abusing notation, we refer to the weak equivalence itself as the model.
The work of Mandell~\cite{mandell:01} shows that if $R$ is the algebraic
closure of $\F{p}$, then the cell model of the cochain algebra of a space
captures the $p$-adic homotopy type of the space.
In particular, the cell model carries all of the information on the
Steenrod operations in cohomology.
This fact, along with the fact that every morphism has a cell model,
addresses the difficulties encountered when trying to use commutative
algebras to model $p$-local maps and spaces.

Let $LX$ and $\Omega X$ be the free and based loop spaces, respectively,
on the pointed, simply-connected topological space $X$.
We prove the following Main Theorem and Important Corollary.

\begin{MT}
Let $m_{X}:(\free{E}{V},d) \xrightarrow{\sim} N^{*}(X)$ be a cell algebra
model.
Then there exists a cell algebra model
$(\free{E}{V \oplus sV},d) \xrightarrow{\sim} N^{*}(LX)$,
and the cell extension
$(\free{E}{V},d)\rightarrow(\free{E}{V \oplus sV},d)$
models the evaluation map
$ev:LX \rightarrow X$.
\end{MT}

\begin{IC}
With the notation of the Main Theorem,
there exists a cell model of the form
$(\free{E}{sV},d) \xrightarrow{\sim} N^{*}(\Omega X)$,
and the quotient map
$(\free{E}{V \oplus sV},d)\rightarrow(\free{E}{sV},d)$
models the inclusion
$\Omega X \subset LX$.
\end{IC}

In fact, the Main Theorem and the Important Corollary hold for any cofibrant
models.
They are stated in full generality as Theorem~\ref{thm:free-loop-model} and
Corollary~\ref{cor:based-loop-model}, respectively.
Similar results have been found by Chataur and
Thomas~\cite{chataur-thomas:02} using an operadic Hochschild complex,
and by Fresse~\cite{fresse:01} using a derived functor of the left
adjoint to a mapping functor.
The nice thing about the Important Corollary, as compared to the usual
situation when working $p$-locally
(see~\cite{baues:98,halperin:92,jrsmith:94}), is that
the construction may be iterated:

\begin{cor}
If $X$ is $q$-connected, then there exists a cell model of the form
$(\free{E}{s^{q}V},d) \xrightarrow{\sim} N^{*}(\Omega^{q}X)$.
\end{cor}

\begin{ex}
A cell algebra model
$\free{E}{V} \xrightarrow{\sim} N^{*}(S^{2n+1})$
determines cell algebra models
$\free{E}{sV} \xrightarrow{\sim} N^{*}(\Omega S^{2n+1})$
and
$\free{E}{s^{2}V} \xrightarrow{\sim} N^{*}(\Omega^{2}S^{2n+1})$.
So all of the Steenrod operations in $H^{*}(\Omega^{2}S^{2n+1};\F{p})$
are lurking somewhere in $\free{E}{V}$.
\end{ex}

The outline of the paper is as follows.
Section~\ref{sec:operads} covers the definitions and basic
properties of operads, algebras, and cell algebras, and
our $E_{\infty}$ operad of choice, the Barratt-Eccles operad.
In Section~\ref{sec:cylinders} we construct explicit cylinder and cone
objects for cell $\op{O}$-algebras for certain operads $\op{O}$ including
the Barratt-Eccles operad.
In Section~\ref{sec:models} we assemble various facts
from~\cite{mandell:01} to show that the pushout of cell models gives a
cell model of the pullback of spaces (Lemma~\ref{lem:mandell}) and use
this result to prove the Main Theorem and hence deduce the Important
Corollary.
We end with some computational examples in Section~\ref{sec:examples}.

The first author thanks the CRM (Barcelona) for its hospitality.
The second author would like to thank the Universit\'e d'Angers
for its hospitality, and the Universidad de M\'alaga for providing
an environment so conducive to improving lemmas.

\section{Operads and their algebras}\label{sec:operads}

In this section we fix notation and recall some definitions.
Throughout the paper wee work over a commutative ground ring $R$.

Differential graded (DG) $R$-modules are $\Z$-graded, and we use the
convention $M^{j}=M_{-j}$.
Define the $R$-free chain complex $I$ by $I_{0} = R\{e',e''\}$,
$I_{1}=R\{s\}$, $ds = e'-e''$.
The augmentation $I \rightarrow R$, sending $e'$ and $e''$ to $1$, is
a quasi-isomorphism.
The \emph{cylinder} on a DG $R$-module $M$ is $IM := I \otimes M$.
Note that $IM = M' \oplus M'' \oplus sM$, where $M' = e'\otimes M \cong M$,
$M'' = e'' \otimes M \cong M$, and $sM = s \otimes M$.
Since $I$ is a semifree $R$-module, the augmentation on $I$ induces a
quasi-isomorphism $IM \xrightarrow{\sim} M$.
The \emph{cone} on $M$ is the quotient $CM := IM / (e'-e'') = M \oplus
sM$, and we have a quasi-isomorphism $CM \xrightarrow{\sim} 0$.
The \emph{suspension} on $M$ is the quotient $sM := CM / M$.
Suspension is an isomorphism of lower degree $+1$: $s:M_{j}\cong
(sM)_{j+1}$.

An \emph{operad} $\op{O}$ consists of DG $R$-modules $\op{O}(n)$,
$n \geq 0$, together with a unit map $R \rightarrow \op{O}(1)$, a right
action of the symmetric group $\Sigma_{n}$ on each $\op{O}(n)$, $n \geq 1$,
and chain maps
\[
    \op{O}(n) \otimes \op{O}(j_{1}) \otimes \cdots \otimes \op{O}(j_{n})
        \rightarrow \op{O}(j_{1} + \cdots j_{n})
\]
for $n\geq 1$, $j_{s}\geq 0$, called the \emph{composition products},
that are required to be associative, unital, and equivariant.
See K\v ri\v z and May~\cite{kriz-may:95} for the details.

An $\op{O}$-\emph{algebra} is a DG $R$-module $A$ together with
chain maps
\[
    \op{O}(n) \otimes A^{\otimes n} \rightarrow A
\]
for $n \geq 0$ that are associative, unital, and equivariant.
Once again, see K\v ri\v z and May~\cite{kriz-may:95}.
Coproducts other than direct sums and tensor products shall be denoted
by the symbol $\vee$.

The \emph{free} $\op{O}$-\emph{algebra} on a DG $R$-module $(V,d)$ is
defined by
\[
    \free{O}{V,d} =
        \bigoplus_{n\geq 0}\left(\op{O}(n) \otimes_{\Sigma_{n}} V^{\otimes n}
            \right).
\]
We will abuse notation and write $(\free{O}{V},d)$ for an $\op{O}$-algebra
that is free if differentials are ignored.

For a simplicial set $Y$, denote by $N_{*}(Y)$ and $N^{*}(Y)$ the normalised
chains and cochains on $Y$, respectively.
If $X$ is a space we take $N_{*}(X) = N_{*}(S(X))$ where $S(X)$ is the
singular simplicial set on $X$.

Let $X$ be a set.
Denote by $W(X)$ the \emph{standard simplicial resolution of} $X$,
where $W(X)_{n} = X^{n+1}$, $n\geq 0$.
Face and degeneracy maps are defined by deletion and repetition, respectively.

The \emph{Barratt-Eccles operad} $\op{E}$ is defined as
$\op{E}(n) = N_{*}(W(\Sigma_{n}))$, with composition products determined
by block permutations.
The reader is referred to Berger and Fresse~\cite{berger-fresse:01} for a
detailed description of the composition product.
For $n \geq 0$, $\op{E}(n)$ is an $R[\Sigma_{n}]$-free resolution of $R$,
so $\op{E}$ is an $E_{\infty}$ operad in the sense of K\v ri\v z and May.
In addition, Berger and Fresse showed that the category of $\op{E}$-algebras
has a particularly nice model structure.
The fibrations are the surjections, the weak equivalences are the
quasi-isomorphisms, and the cofibrations are retracts of cell extensions,
which we now define.

Let $\op{O}$ be an operad.
An $\op{O}$-algebra morphism $j:A \rightarrow B$ is called a
\emph{cell extension} (\emph{relative cell inclusion} in the language
of Mandell~\cite{mandell:01}) if
\begin{enumerate}
\item forgetting differentials, $B \cong A \vee \free{O}{V}$,
\item $j$ is the canonical inclusion, and
\item $V$ is the union of a nested sequence of submodules $V(k)$,
$k \geq 0$, such that
$V(0)$ and $V(k)/V(k-1)$, $k \geq 0$, are $R$-free, $dV(0) \subset A$,
and $d(V(k)) \subset A \vee \free{O}{V(k-1)}$ for $k \geq 1$.
\end{enumerate}
In particular, we may write $V(k) = V_{k} \oplus V(k-1)$ with
$V_{k} \cong V(k)/V(k-1)$ and $d(V_{k}) \subset A \vee \free{O}{V(k-1)}$.
A \emph{cell algebra} is a cell extension of $R$.
A \emph{cell model} of a morphism $\varphi:A \rightarrow A'$ is a
factorisation of $\varphi$ as a cell extension
followed by a weak equivalence.
A cell model of an $\op{O}$-algebra $A$ is a cell model of the unit morphism
$R \rightarrow A$.

Let $A$ be an element of a closed model category, and fix a cylinder object
$IA$.
Set $LA = A \vee_{A \vee A} IA$ and
$SA = R \vee_{A} LA$, where $R$ is the terminal object in the
category.
Both $LA$ and $SA$ depend upon choice of cylinder object, but only
up to weak equivalence.

\section{Cylinder objects and acyclic closures}\label{sec:cylinders}

In this section we construct an explicit cylinder object and acyclic closure
for a given cell algebra.
Let $\op{O}$ be an operad such that the $\op{O}$-algebras form a closed
model category where the fibrations are the surjections and the weak
equivalences are the quasi-isomorphisms.
This condition is satisfied, for example, by the associative algebra
operad, the commutative algebra operad if $R \supseteq \Q$, the
Barratt-Eccles operad, or any cofibrant operad.
Essentially, we want the free $\op{O}$-algebra functor to preserve
quasi-isomorphisms.

\begin{prop}\label{prop:cylinder}
Let $A = (\free{O}{V},d)$ be a cell $\op{O}$-algebra.
Set $V' = V'' = V$.
Then the fold map $\nabla:A \vee A \rightarrow A$ has a surjective cell
model
\[
  \xymatrix{
    A \vee A \ar@{>->}[r]
    & IA \ar@{->>}[r]^{\sim}
    & A
  }
\]
where
$IA = (\free{O}{V' \oplus V'' \oplus sV},d)$.
Furthermore, if $x \in V_{k}$, then
\[
  d(sx) - x' + x'' \in \free{O}{V'(k-1) \oplus V''(k-1) \oplus sV(k-1)}.
\]
\end{prop}

\begin{proof}
We construct $IA$ recursively, proceeding along the filtration on $V$.
First we introduce some notation.  For $k \geq 0$, let
$V'(k)$ and $V''(k)$ be isomorphic copies of $V(k)$.
The images of $x \in V(k)$ in $V'(k)$ and $V''(k)$ will be denoted
$x'$ and $x''$, respectively.
Set $A(k)=(\free{O}{V(k)},d)$ and
$IA(k) = (\free{O}{V'(k) \oplus V''(k) \oplus sV(k)},d)$.
Write $A'(k)$ and $A''(k)$ for the subalgebras of $IA(k)$ generated
by $V'(k)$ and $V''(k)$, respectively.

Suppose that we have constructed the algebra $IA(k)$ such that the
epimorphism
$\eta_{k}:IA(k)\rightarrow A(k)$ defined by
$\eta_{k}(x')=\eta_{k}(x'')=x$,
$\eta_{k}(sx)=0$ for all $x \in V(k)$, is a weak equivalence.
Recall that $V(k+1)=V_{k+1}\oplus V(k)$, with
$d(V_{k+1}) \subset A(k)$.
Let $x$ be a basis element of $V_{k+1}$.
By definition, $dx'\in A'(k)$, $dx'' \in A''(k)$, and
$\eta_{k}(d(x'-x''))=0$.
Since $\eta_{k}$ is a trivial fibration, $\ker{\eta_{k}}$ is
contractible,
and so the cycle $d(x'-x'')=d\Phi$ for some
$\Phi \in \ker{\eta_{k}} \subset IA(k)$.
Extend the differential to $sV(k+1)$ by setting
$dsx = x' - x'' - \Phi$.
Set $\eta_{k+1}(v') = \eta_{k+1}(v'') = v$ and $\eta_{k+1}(sv) = 0$
to extend $\eta_{k}$ to $\eta_{k+1}$ .

We need to show that $\eta_{k+1}$ is a weak equivalence.
To this end, define filtrations on $IA(k+1)$ and $A(k+1)$ by
$F_{0} = IA(k+1)$,
$F_{1} = IA(k)$,
and $G_{0} = A(k+1)$, $G_{1} = A(k)$.
The morphism $\eta_{k+1}$ preserves the filtrations and so defines a morphism
of strongly convergent right-half-plane spectral sequences.
At the $E_{0}$-term, the induced morphism is
\[
  E_{0}(\eta_{k+1})=\free{O}{\epsilon}\vee\eta_{k}:
    \free{O}{IV_{k}} \vee IA(k) \rightarrow
    \free{O}{V_{k}} \vee A(k)
\]
where $IV_{k}$ is the cylinder on $(V_{k},0)$.
Since $IV_{k} \xrightarrow{\sim} V_{k}$ and $\free{O}{-}$ preserves weak
equivalences, $E_{0}(\eta_{k+1})$ is a weak equivalence.
\end{proof}

The \emph{acyclic closure}, or \emph{cone}, of $A$ is the algebra
$CA := R \vee_{A} IA$.

\begin{cor}\label{cor:acyclic-closure}
$A \rightarrow CA$ is a cell extension, and the augmentation
$A \rightarrow R$ factors as
$A \rightarrowtail CA \stackrel{\sim}{\twoheadrightarrow} R$.
\end{cor}

\begin{proof}
Apply $R \vee_{A}-$ to
$A \vee A \rightarrowtail IA \stackrel{\sim}{\twoheadrightarrow} A$.
\end{proof}

\section{Models for free and based loop spaces}\label{sec:models}

In this section we specialise to an $E_{\infty}$-operad $\op{E}$,
such as the Barratt-Eccles operad.
We assemble facts from Mandell~\cite{mandell:01} to prove
Lemma~\ref{lem:mandell}.
We then use a characterisation of $LX$ as a pullback to prove
Theorem~\ref{thm:free-loop-model} and deduce
Corollary~\ref{cor:based-loop-model}.

Let $MX$ be the space of free paths on the simply-connected space $X$.
We can describe $LX$ as the pullback of the diagram
\[
    \xymatrix{
        LX \ar[r] \ar[d]^{ev} \ar[r]    & MX \ar[d]^{(e_{0},e_{1})} \\
        X  \ar[r]^(0.4){\Delta} & X \times X
    }
\]
where $ev$ is evaluation at the basepoint of the loop and $e_{0}$,
$e_{1}$ are evaluations at the endpoints of the path.

\begin{thm}\label{thm:free-loop-model}
Let $m_{X}:A \xrightarrow{\sim} N^{*}(X)$ be a cofibrant model.
Then the composite
\[
  A \xrightarrow{\sim} N^{*}(X) \xrightarrow{ev^{*}} N^{*}(LX)
\]
has a cofibrant model
\[
  A \rightarrowtail LA \stackrel{\sim}{\rightarrow} N^{*}(LX)
\]
where $LA := A \vee_{A \vee A}IA$.
If $A = (\free{E}{V},d)$ is a cell algebra, then
$LA = (\free{E}{V \oplus sV},d)$ and $A \rightarrowtail LA$ is a
cell extension.
\end{thm}

As mentioned above, the following lemma is implicit in
Mandell~\cite{mandell:01}.

\begin{lem}\label{lem:mandell}
Let $\pi:E \rightarrow B$ be a fibration, and let
$f:X \rightarrow B$ be continuous.
Suppose there exists a commutative diagram
\[
  \xymatrix{
    A_{X} \ar[d]^{\sim}_{\theta_{X}}
      & A_{B} \ar[d]^{\sim}_{\theta_{B}} \ar[l]_{\theta_{f}} \ar@{>->}[r]
      & A_{E} \ar[d]^{\sim}_{\theta_{E}} \\
    N^{*}(X)
      & N^{*}(B) \ar[l] \ar[r]
      & N^{*}(E)
  }
\]
of $E_{\infty}$ algebras and morphisms, in which the algebras in the top
row are cofibrant.
Then the map induced by pushout
\[
  \theta:A_{X} \vee_{A_{B}} A_{E} \rightarrow N^{*}(X \times_{B} E)
\]
is a weak equivalence.
\end{lem}

\begin{proof}
First we suppose that $\theta_{f}$ is a cofibration, and that all
vertical morphisms are fibrations.
Following Mandell~\cite{mandell:01}, let
$N(\beta(A_{X},A_{B},A_{E}))$ denote the normalised chains on the
simplicial bar construction.
The composition of natural maps
\[
  N(\beta(A_{X},A_{B},A_{E})) \rightarrow A_{X} \vee_{A_{B}} A_{E}
    \xrightarrow{\theta} N^{*}(X \times_{B} E)
\]
is a weak equivalence by~\cite[Lemma 5.2]{mandell:01}, while the first
map is a weak equivalence by~\cite[Theorem 3.5]{mandell:01}.
By the two-out-of-three rule, $\theta$ is a weak equivalence.

If the vertical arrows are not necessarily fibrations, use the closed
model category structure on $E_{\infty}$ algebras to form the diagram
\[
  \xymatrix{
    A_{X} \ar[dd]^{\sim}_{\theta_{X}}
      & & A_{B} \ar[d]^{\sim}_{\psi_{B}} \ar[ll]_{\theta_{f}}
        \ar@{>->}[r]
      & A_{E} \ar[d]^{\sim}_{\psi_{E}} \\
    & B_{X} \ar@{->>}[dl]^{\sim}_{\varphi_{X}}
      & B_{B} \ar@{->>}[d]^{\sim}_{\varphi_{B}}
          \ar@{>->}[l]_{\varphi_{f}} \ar@{>->}[r]^{\varphi_{\pi}}
      & B_{E} \ar@{->>}[d]^{\sim}_{\varphi_{E}} \\
    N^{*}(X)
      & & N^{*}(B) \ar[ll]_{f^{*}} \ar[r]^{pi^{*}}
      & N^{*}(E)
  }
\]
in which $\varphi_{B}\circ\psi_{B}$,
$\varphi_{X}\circ\varphi_{f}$ and
$\varphi_{E}\circ\varphi_{\pi}$ are factorisations of
$\theta_{B}$, $f^{*}\circ\varphi_{B}$ and $\pi^{*}\circ\varphi_{B}$,
respectively,  into a cofibration followed by a trivial fibration,
and $\psi_{E}$ is a lift for the diagram
\[
  \xymatrix{
    A_{B} \ar[r] \ar@{>->}[d] & B_{E} \ar@{->>}[d]^{\sim} \\
    A_{E} \ar[r] & N^{*}(E).
  }
\]
By the first paragraph of the proof, the pushout morphism
$\varphi : B_{X} \vee_{B_{B}} B_{E}
    \rightarrow N^{*}(X \times_{B} E)$
is a weak equivalence.

If $\theta_{X}$ is a surjection, then factor $\varphi_{f}\circ\psi_{B}$
as $\xi_{X}\circ\psi_{f}$, where $\psi_{f}:A_{B} \rightarrowtail C_{X}$
is a cofibration and
$\xi_{X}:C_{X} \stackrel{\sim}{\twoheadrightarrow} B_{X}$ is a trivial
fibration.
The pushout morphism
$\tilde{\psi} : C_{X} \vee_{A_{B}} A_{E}
    \rightarrow B_{X} \vee_{B_{B}} B_{E}$
induced by $\xi_{X}$, $\psi_{B}$, and $\psi_{E}$, is a weak equivalence
by~\cite[Theorem 3.2]{mandell:01}.
Now lift $\varphi_{X}\circ\xi_{X}$ through $\theta_{X}$ to define a
weak equivalence
$\eta_{X}:C_{X}\xrightarrow{\sim} A_{X}$ such that
$\theta_{f}=\eta_{X}\circ\psi_{f}$.
The pushout morphism
$\tilde{\eta} : C_{X} \vee_{A_{B}} A_{E}
    \rightarrow A_{X} \vee_{A_{B}} A_{E}$
induced by $\eta_{X}$ is a weak equivalence
by~\cite[Theorem 3.2]{mandell:01}.
By uniqueness of pushout,
$\theta\circ\tilde{\eta}=\varphi\circ\tilde{\psi}$.
It follows that $\theta$ is a weak equivalence.

If $\theta_{X}$ is not necessarily surjective, then factor it as
$\theta_{X} = p_{X}\circ i_{X}$, where
$p_{X}:U_{X}\rightarrow A_{X}$ is a trivial fibration
and $i_{X}$ is a cofibration.  By the two-out-of-three rule, $i_{X}$ is
a weak equivalence.  By the previous paragraph, $p_{X}$, $\theta_{B}$,
and $\theta_{E}$ induce a weak equivalence
$\tilde{\theta}:U_{X}\vee_{A_{B}} A_{E}
    \rightarrow N^{*}(X \times_{B} E)$.
Furthermore, the map of pushouts
$i: A_{X}\vee_{A_{B}} A_{E} \rightarrow
  U_{X}\vee_{A_{B}} A_{E}$
induced by $i_{X}$ and the identities on $A_{B}$ and $A_{E}$ is a weak
equivalence  by~\cite[Theorem 3.2]{mandell:01}.
Thus $\theta = \tilde{\theta}\circ i$ is a weak equivalence.
\end{proof}

Setting $X = *$ we obtain the $E_{\infty}$ analogue of the theorem on
the model of the 
fibre~\cite{grivel:79,halperin:83,dupont-hess:94,ndombol:98,menichi2:01}.

\begin{cor}\label{cor:fibre-model}(cf. Chataur~\cite{chataur:02})
Let $\pi:E\rightarrow B$ be a fibration with fibre $F$.
Then $R \vee_{A_{B}} A_{E}$ is a model for $N^{*}(F)$.
\end{cor}

Of course, Corollary~\ref{cor:fibre-model} applies to the homotopy fibre
of an arbitrary map.

\begin{lem}\label{lem:product-model}
Let $m_{X}:A \xrightarrow{\sim} N^{*}(X)$ be a cofibrant model.
Then the fold map $A \vee A \xrightarrow{\nabla} A$ models
$\Delta^{*}$.
That is, there exists a weak equivalence
$m_{X \times X}:A \vee A \xrightarrow{\sim} N^{*}(X \times X)$ such
that the diagram
\[
    \xymatrix{
        A \vee A \ar[d]^{\sim} \ar[r]^{\nabla}  & A \ar[d]^{\sim}   \\
        N^{*}(X \times X) \ar[r]^{\Delta^{*}}   & N^{*}(X)
    }
\]
commutes.
\end{lem}

\begin{proof}
Denote by $\pi_{1},\pi_{2}:X \times X \rightarrow X$ the canonical
projections.

The product $X \times X$ is the pullback of the diagram
$X \rightarrow * \leftarrow X$.
Since $N^{*}(*)=R$ and $m_{X}$ commutes with unit morphisms,
Lemma~\ref{lem:mandell} states that the pushout morphism
$m_{X \times X}:A_{X} \vee A_{X}
    \xrightarrow{\sim} N^{*}(X \times X)$ is a weak equivalence.
Since $\Delta\circ\pi_{i} = 1_{X}$ for $i = 1,2$,
the compositions
$m_{X}\circ\nabla$ and
$\Delta^{*}\circ m_{X \times X}$ both make the diagram
\[
    \xymatrix{
        R \ar@{>->}[r] \ar@{>->}[d] & A_{X} \ar@{->>}[r]^{\sim}
            \ar[d]  & N^{*}(X) \ar[d]^{\pi_{1}^{*}} \\
        A_{X} \ar[r] \ar@{->>}[d]^{\sim}
            & A_{X} \vee A_{X} \ar@{.>}[dr]
            & N^{*}(X \times X) \ar[d]^{\Delta^{*}} \\
        N^{*}(X) \ar[r]^{\pi_{2}^{*}}
            & N^{*}(X \times X) \ar[r]^{\Delta^{*}}
            & N^{*}(X)
    }
\]
commute.  Therefore
$m_{X}\circ\nabla=\Delta^{*}\circ m_{X \times X}$.
\end{proof}

If $A$ is an $E_{\infty}$ algebra, denote by $IA$ a cylinder object on
$A$, and $j_{1},j_{2}:A \rightarrow IA$ the canonical inclusions.

\begin{lem}\label{lem:path-space-model}
Let $e = (e_{0},e_{1}):MX \rightarrow X \times X$.
Given the hypotheses of Lemma~\ref{lem:product-model},
the canonical morphism $A \vee A \xrightarrow{j_{1}+j_{2}} IA$ is a
cofibrant model for $e^{*}$.  That is, there exists a weak
equivalence $m_{MX}:IA \xrightarrow{\sim} N^{*}(MX)$ such that the
diagram
\[
    \xymatrix{
        A \vee A \ar[d]^{\sim}_{m_{X \times X}} \ar@{>->}[r]^{j_{1}+j_{2}}
    & IA \ar[d]^{\sim}_{m_{MX}} \\
        N^{*}(X \times X) \ar[r]^{e^{*}}    & N^{*}(MX)
    }
\]
commutes.
If $A$ is a cell algebra, then $IA$ may be chosen to be the cylinder of
Proposition~\ref{prop:cylinder}.
\end{lem}

\begin{proof}
To construct $m_{MX}$, we note that $e_{0}\simeq e_{1}:MX \rightarrow X$.

Therefore $e_{0}^{*}\circ m_{X}\simeq e_{1}^{*}\circ m_{X}$ as
$E_{\infty}$ morphisms, and so the composite
\[
  A \vee A \xrightarrow{m_{X} \vee m_{X}}
    N^{*}(X) \vee N^{*}(X)
    \xrightarrow{e_{0}^{*} + e_{1}^{*}} N^{*}(MX)
\]
factors as
$A \vee A \xrightarrow{j_{1}+j_{2}} IA
  \xrightarrow{m_{MX}} N^{*}(MX)$.

Since $e_{0}$ is a homotopy equivalence,
$e_{0}^{*}:N^{*}(X)\rightarrow N^{*}(MX)$ is a weak equivalence.
Then $m_{MX}$ is a weak equivalence because
$m_{MX}\circ j_{1} = e_{0}^{*} \circ m_{X}$.

Now,
$e_{0}^{*} + e_{1}^{*} = e^{*} \circ (\pi_{1}^{*} + \pi_{2}^{*})$
because
$e_{i} = \pi_{i+1}\circ(e_{0},e_{1})$ for $i = 0,1$,
and
$(\pi_{1}^{*} + \pi_{2}^{*}) \circ (m_{X}\vee m_{X}) = m_{X \times X}$
by uniqueness of pushout.
We deduce that
$e^{*} \circ m_{X \times X}
  = (e_{0}^{*} + e_{1}^{*}) \circ (m_{X}\vee m_{X})$.

The above proof works for any cylinder object on $A$, in particular the
cylinder of Proposition~\ref{prop:cylinder} if $A$ happens to be
a cell algebra.
\end{proof}

Now it is time for the

\begin{proof}[Proof of Theorem~\ref{thm:free-loop-model}]
Consider the diagram of continuous maps
\begin{equation}\label{eq:LX}
    \xymatrix{
        X \ar[r]^{\Delta}
            & X \times X
            & MX \ar[l]_{e}
    }
\end{equation}
where $e = (e_{0},e_{1})$ and $MX$ is the free path space on $X$.
Recall that the pullback of (\ref{eq:LX}) is the free loop space $LX$.
Use Lemmas~\ref{lem:product-model} and~\ref{lem:path-space-model}
to form the diagram
\[
    \xymatrix{
  A \ar@{->>}[d]^{m_{X}}_{\sim}
    & A \vee A \ar[l]_{\nabla} \ar@{>->}[r]
      \ar[d]^{m_{X \times X}}_{\sim}
    & IA \ar[d]^{m_{MX}}_{\sim} \\
        N^{*}(X)
            & N^{*}(X \times X) \ar[l]_{\Delta^{*}}\ar[r]^{e^{*}}
            & N^{*}(MX).
    }
\]
By Lemma~\ref{lem:mandell}, the pushout morphism
$m: A \vee_{A \vee A} IA
  \rightarrow N^{*}(LX)$ is a weak equivalence.
Furthermore, if
$j_{A}:A \rightarrowtail A \vee_{A \vee A} IA$
is the natural map, then
$m\circ j_{B} = ev^{*} \circ m_{X}$ by construction.

The natural map $A \rightarrow LA$ is a cofibration since
$A \vee A \rightarrowtail IA$ is one.
If $A = (\free{E}{V},d)$ is a cell algebra, with cellular filtration
$V(k)$, then
the filtration $(sV)(k) = s(V(k))$ exhibits $A \rightarrow LA$ as
a cell extension, completing the proof.
\end{proof}

\begin{cor}\label{cor:based-loop-model}
There exists a cell model of the form
$(\free{E}{sV},d) \xrightarrow{\sim} N^{*}(\Omega X)$.
\end{cor}

\begin{proof}[Proof of Corollary~\ref{cor:based-loop-model}]
The based loop space $\Omega X$ is the fibre of the evaluation map
$ev:LX \rightarrow X$.
By Corollary~\ref{cor:fibre-model}, the induced map
$m_{\Omega X}:R \vee_{A} LA \rightarrow N^{*}(\Omega X)$ is a weak
equivalence.
If $A = \free{E}{V}$ is a cell algebra, then we may assume
$LA = \free{E}{V \oplus sV}$, is a cell extension of $A$, whence
$SA := R \vee_{A} LA = \free{E}{sV}$ is a cell algebra.
\end{proof}

\begin{rmk}
We may choose to characterise $\Omega X$ as the fibre of two other fibrations,
namely, the evaluation maps $e:MX \rightarrow X \times X$
and $e_{0}:PX \rightarrow X$, where $PX$ is the space of paths ending at
$x_{0}$.
Of course, $ev$ is the pullback of $e$ along the diagonal
$\Delta:X \rightarrow X \times X$ and $e_{0}$ is the pullback of $e$
along $(id,x_{0}):X \rightarrow X \times X$.
The resulting models form the diagram of pushout squares
\[
    \xymatrix{
    A \ar@{>->}[d]
        & A \vee A \ar[l] \ar[r] \ar@{>->}[d]
        & A \ar@{>->}[d] \\
    CA
        & IA \ar[l] \ar[r]
        & LA
    }
\]
shows that $SA = R \vee_{A \vee A}IA = R \vee_{A} CA$.
In practice, $R \vee_{A} CA$ is the easiest to compute.
\end{rmk}

\section{Examples}\label{sec:examples}

In all the examples below we work over $\F{2}$. Recall that
$\op{E}(2)_{q} = \F{2}\pi\cdot e_{q}$, where $\pi$ is the cyclic
group of order $2$ with generator $\tau$, and $de_{q} = (1 +
\tau)e_{q-1}$.

\subsection{Free and based loops on an Eilenberg-Mac Lane space}\label{ex:free-based-loops}
We calculate the models for the free and based loop spaces on
$K(\Zmod{2},n)$. We show that the model for the based loop space
is the model for $K(\Zmod{2},n-1)$, and that the model for the
free loop space splits as expected.

A model for $N^{*}(K(\Zmod{2},n);\F{2})$ is
$A(n)=(\free{E}{u,v},dv=u+e_{n}(u,u))$ (see~\cite{mandell:01}).
Then
$IA(n) = \free{E}{u',v',u'',v'',su,sv}$,
with
$dsu=u'+u''$.
Since
\[
    d(v'+v'')=d(su + e_{n}(u'',su) + e_{n}(su,u') + e_{n-1}(su,su)),
\]
we may take
\[
    dsv = v' + v'' + su + e_{n}(u'',su) + e_{n}(su,u') + e_{n-1}(su,su).
\]
Setting
$u'=u''$
and
$v'=v''$,
we find
$LA(n) = \free{E}{u,v,su,sv}$,
with
$dsu=0$
and
$dsv = su + e_{n}(u,su) + e_{n}(su,u) + e_{n-1}(su,su)$.
Then
\[
    SA(n) = (\free{E}{su,sv},dsv = su + e_{n-1}(su,su))
\]
is evidently $A(n-1)$, a model for
$N^{*}(K(\Zmod{2},n-1);\F{2})$.
The natural map $LA(n) \rightarrow SA(n)$
has a splitting, defined by $su \mapsto su$,
$sv \mapsto sv + e_{n+1}(u,su)$,
that defines an isomorphism
$A(n) \vee SA(n) \xrightarrow{\cong} LA(n)$.

\subsection{The algebraic model of an elementary fibration}\label{ex:elem-fib}
Let $X$ be a ``nice'' topological space, i.e. it is $2$-complete
and of finite $2$-type.
We want to completely determine the
algebraic model of an elementary fibration over $X$. That is to
say, we consider an algebraic model of the pullback square of the
path fibration over $K(\Zmod{2},n)$ over a map
$f:X\rightarrow K(\Zmod{2},n)$,
\[
    \xymatrix{
        X' \ar[r] \ar[d] & PK(\Zmod{2},n) \ar[d] \\
        X \ar[r]^{f} & K(\Zmod{2},n)
    }
\]
By Lemma~\ref{lem:mandell}, the model is given by the pushout square of
algebras:
\begin{equation}\label{eq:elem-fib}
    \xymatrix{
        A_{X'} & CA(n) \ar[l] \\
        A_{X} \ar[u] & A(n). \ar[l] \ar[u]
    }
\end{equation}
where $A(n)$ is the model for $N^{*}(K(\Zmod{2},n);\F{2})$ given in
Example~\ref{ex:free-based-loops}, and $CA(n)$ is its acyclic closure
(recall Corollary~\ref{cor:acyclic-closure}).
The acyclic closure is a cell extension that models the path fibration
\[
    \xymatrix{
        K(\Zmod{2},n-1)\ar[r]
            & PK(\Zmod{2},n)\ar[r] & K(\Zmod{2},n).}
\]
and its cofibre $SA(n)$ is $A(n-1)$, a model for
$N^{*}(K(\Zmod{2},n-1);\F{2})$.
The differential in $CA(n)$ extends that of $A(n)$, with
$dsu=u$ and
$dsv=v+su+e_{n-1}(su,su)+e_n(su,u)$.

We can now give a model of the map $X'\rightarrow X$.
Let $\phi:A(n)\rightarrow A_X$ be a model for $f$.
Set $z=\phi(u)$ and $z'=\phi(v)$.
By the diagram (\ref{eq:elem-fib}),
$A_{X'}=A_X \vee A(n-1) = A_{X} \vee \free{E}{su,sv}$,
and the differential in $A_{X'}$ satisfies
$dsu=z$
and
\[
    dsv=u+e_{n-1}(su,su)+e_n(su,z)+z'.
\]

\subsection{Loop spaces of some Postnikov $2$-towers.}
We are interested by determining a model of the loop space
of the space $X_1$ obtained as the total space of the following elementary
fibration:
\[
    \xymatrix{
 X_1 \ar[r] \ar[d]
        & PK(\Zmod{2},n+p) \ar[d] \\
 K(\Zmod{2},n) \ar[r]^{f}
        & K(\Zmod{2},n+p)
    }
\]
where $f$
represents a Steenrod square $\Sq{p}$.
Recall the algebra $A(k) = \free{E}{u_{k},v_{k-1}}$ defined in
Example~\ref{ex:free-based-loops}.
A model for $f$ is given by
\[
    \phi:A(n+p)\rightarrow A(n)
\]
where
$\phi(u_{n+p})=e_{n-p}(u_n,u_n)$ and
$\phi(v_{n+p-1})=\gamma_{n+p-1}$.
Here $\gamma_{n+p-1}$ is an element that satisfies
$d\gamma_{n+p-1}=e_{n-p}(u_n,u_n)+e_{n+p}(e_{n-p}(u_n,u_n),e_{n-p}(u_n,u_n))$.
Applying the formulas of the preceding examples we obtain
\[
    A_{X_1}=\free{E}{u_n,v_{n-1},w_{n+p-1},t_{n+p-2}}
\]
together with the differential
\begin{eqnarray*}
    du_n        & = & 0 \\
    dv_{n-1}    & = & u_n+e_n(u_n,u_n) \\
    dw_{n+p-1}  & = & e_{n-p}(u_n,u_n) \\
    dt_{n+p-2}  & = & w_{n+p-1}+\gamma_{n+p-1}+e_{n+p-1}(w_{n+p-1},w_{n+p-1})\\
                &   & +e_{n+p}(e_{n-p}(u_n,u_n),w_{n+p-1}).
\end{eqnarray*}
A differential on the acyclic closure of $A_{X_1}$,
\[
    CA_{X_1}=
        \free{E}{u_n,v_{n-1},w_{n+p-1},t_{n+p-2},
            u'_{n-1},v'_{n-2},w'_{n+p-2},t'_{n+p-3}}
\]
is given by
\begin{eqnarray*}
    du'_{n-1}   & = & u_n \\
    dv'_{n-2}   & = & v_{n-1}+u'_{n-1}+e_{n-1}(u'_{n-1},u'_{n-1})
                 +e_{n}(u_n,u'_{n-1}) \\
    dw'_{n+p-2} & = & w_{n+p-1}+e_{n-p}(u'_{n-1},u_n)
                    +e_{n-p-1}(u'_{n-1},u'_{n-1}) \\
    dt'_{n+p-3} & = & t_{n+p-2}+w'_{n+p-2}+e_{n+p-2}(w'_{n+p-2},w'_{n+p-2}) \\
                &   & + e_{n+p-1}(e_{n-p-1}(u'_{n-1},u'_{n-1}),w'_{n+p-2})+
                    \gamma'_{n+p-2}
\end{eqnarray*}
such that in the cofibre
$SA_{X_{1}}=\free{E}{u'_{n-1},v'_{n-2},w'_{n+p-2},t'_{n+p-3}}$
we have:
\[
    d\gamma'_{n+p-2}
        =e_{n-p-1}(u'_{n-1},u'_{n-1})
        +e_{n+p-1}(e_{n-p-1}(u'_{n-1},u'_{n-1}),e_{n-p-1}(u'_{n-1},u_{n-1})).
\]
Hence $SA_{X_{1}}$ models $N^{*}(\Omega X_{1};\F{2})$.
The differential is
\begin{eqnarray*}
    du'_{n-1}   & = & 0, \\
    dv'_{n-2}   & = & u'_{n-1}+e_{n-1}(u'_{n-1},u'_{n-1}), \\
    dw'_{n+p-2} & = & e_{n-p-1}(u'_{n-1},u'_{n-1}), \\
    dt'_{n+p-3} & = & w'_{n+p-2}+e_{n+p-2}(w'_{n+p-2},w'_{n+p-2}) \\
                &   & + e_{n+p-1}(e_{n-p-1}(u'_{n-1},u'_{n-1}),w'_{n+p-2})+
                    \gamma'_{n+p-2}.
\end{eqnarray*}
Note that we have recovered algebraically the fact that $\Omega X_1$ can be
obtained as the pull-back
\[
    \xymatrix{
 \Omega X_1 \ar[r] \ar[d]
        & PK(\Zmod{2},n+p-1) \ar[d] \\
 K(\Zmod{2},n-1) \ar[r]
        & K(\Zmod{2},n+p-1)
    }
\]
where the bottom horizontal map is also $\Sq{p}$.

\bibliographystyle{amsplain}
\bibliography{homotopy}

\providecommand{\bysame}{\leavevmode\hbox to3em{\hrulefill}\thinspace}
\providecommand{\MR}{\relax\ifhmode\unskip\space\fi MR }
\providecommand{\MRhref}[2]{%
  \href{http://www.ams.org/mathscinet-getitem?mr=#1}{#2}
}
\providecommand{\href}[2]{#2}

\end{document}